\documentclass{llncs}

\def\Z{{\rm \mathbf{Z}}}
\def\C{{\rm \mathbf{C}}}
\def\R{{\rm \mathbf{R}}}
\def\Q{{\rm \mathbf{Q}}}
\def\N{{\rm \mathbf{N}}}

\def\e{{\rm e}}

\begin{document}

\pagestyle{empty}

\mainmatter
\title{A discrete Farkas lemma}
\titlerunning{A discrete Farkas lemma}  
\author{Jean B. Lasserre}
\authorrunning{Jean B. Lasserre}   % abbreviated author list (for running head)
\institute{LAAS-CNRS\\ 7 Avenue du Colonel Roche, 31077 Toulouse cedex
4, France.\\
\email{lasserre@laas.fr}\\
\texttt{http://www.laas.fr/\homedir lasserre}}

\maketitle

\begin{abstract}
Given $A\in \Z^{m\times n}$ and $b\in\Z^m$, we consider the 
issue of existence of a nonnegative integral solution 
$x\in \N^n$ to the system of linear equations
$Ax=b$. We provide a discrete and explicit analogue
of the celebrated Farkas lemma for linear systems in $\R^n$
and prove that checking existence of integral solutions reduces to solving
an explicit linear programming problem of fixed dimension, known in advance. 
\end{abstract}

\section{Introduction}
Let $A\in\Z^{m\times n},b\in\Z^m$ and consider the problem of existence of 
a solution $x\in\N^n$ of the system of linear equations
\begin{equation}
\label{intro1}
Ax\,=\,b,
\end{equation}
that is, the existence of {a \it nonnegative integral} solution of the linear system $Ax=b$.
\vspace{0.2cm}

{\bf Contribution}.
The celebrated {\it Farkas Lemma} in linear algebra states that
\begin{equation}
\label{farkas}
\{x\in\R^n_+\,|\,Ax=b\}\,\neq\,\emptyset\,\Leftrightarrow\,
\left[\,u\in\R^m \mbox{ and }A'u\,\geq\,0\,\right]
\,\Rightarrow\,b'u\,\geq\,0
\end{equation}
(where $A'$ (resp. $b'$) stands for the transpose of $A$ (resp. $b$)).

To the best of our knowledge, there is no {\it explicit} discrete
analogue of (\ref{farkas}). Indeed, the (test) Gomory and Chv\'atal functions
used by Blair and Jeroslow in \cite{blair}
(see also Schrijver in \cite[Corollary 23.4b]{schrijver}) are defined implicitly
and recursively, and do not provide a test directly in terms of the data $A,b$.

In this paper we provide a {\it discrete} and {\it explicit} analogue of Farkas Lemma
for (\ref{intro1}) to have a solution $x\in\N^n$. 
Namely, when $A$ and $b$ have nonnegative entries,
that is, when $A\in\N^{m\times n},b\in\N^m$, we prove that (\ref{intro1})
has a solution $x\in\N^n$ {\it if and only if} the polynomial
$z\mapsto z^b-1$ ($:=z_1^{b_1}\cdots z_m^{b_m}-1$) of $\R[z_1,\ldots,z_m]$, can be written
\begin{equation}
\label{intro2}
z^b-1\,=\,\sum_{j=1}^nQ_j(z)(z^{A_j}-1)\,=\,
\sum_{j=1}^nQ_j(z)(z_1^{A_{1j}}\cdots z_m^{A_{mj}}-1)
\end{equation}
for some polynomials $\{Q_j\}$ in $\R[z_1,\ldots,z_m]$ with {\it nonnegative} coefficients.
In other words,
\begin{equation}
\label{newfarkas}
\{x\in\N^n\,|\,Ax=b\}\,\neq\,\emptyset\,\Leftrightarrow\,
z^b-1\,=\,\sum_{j=1}^nQ_j(z)(z^{A_j}-1),
\end{equation}
for some polynomials $\{Q_j\}$ in $\R[z_1,\ldots,z_m]$ with {\it nonnegative} coefficients.
(Of course, the {\it if} part of the equivalence in (\ref{newfarkas})
is the hard part of the proof.)

Moreover, the degree of the $Q_j$'s
is bounded by $b^*:=\displaystyle{\sum_{j=1}^mb_j-\min_k\sum_{j=1}^mA_{jk}}$.

Therefore, checking the existence of a solution $x\in\N^n$ to
$Ax=b$, reduces to checking whether or not  there is a 
{\it nonnegative} solution $y$ to a system
of {\it linear} equations where (i) $y$ is the vector of unknown
nonnegative coefficients 
of the $Q_j$'s and (ii), the (finitely many) linear equations
identify coefficients of same power in both sides of (\ref{intro2}).
This is a linear programming (LP) problem with
$n s(b^*)$ variables and $s(b^*+\max_k\sum_jA_{jk})$ constraints, where
$s(u):={m+u\choose u}$ denotes the 
dimension of the vector space of polynomials of degree $u$ in $m$ variables.
In addition, all the coefficients of the associated matrix
of constraints are all $0$ or $\pm 1$. 
For instance, checking the existence of a solution $x\in\N^n$ to the knapsack equation
$a'x=b$, reduces to solving a LP problem with $n(b+1-\min_ja_j)$ variables and $b+1
+\max_ja_j-\min_ja_j$ equality constraints. 
This result is also extended to the case where $A\in\Z^{m\times n},b\in\Z^m$, that is, when $A$
and $b$ may have nonnegative entries.

We call (\ref{newfarkas}) a {\it Farkas lemma} because as (\ref{farkas}), it states a condition
on the {\it dual} variables $z$ associated with the constraints $Ax=b$.
In addition, let $z:=\e ^\lambda$ and notice that the
basic ingredients $b'\lambda$ and $A'\lambda$ of (\ref{farkas}), also 
appear in (\ref{newfarkas}) via $z^b$ which becomes $e^{b'\lambda}$ and
via $z^{A_j}$ which becomes $\e^{(A'\lambda)_j}$.
Moreover, if indeed $z^b-1$ has the representation (\ref{newfarkas}),
then  whenever $\lambda\in\R^m$ and $A'\lambda\geq0$ (letting $z:=\e^\lambda$)
\[\e^{b'\lambda}-1=\sum_{j=1}^nQ_j(\e^{\lambda_1},\ldots,e^{\lambda_m})
\left[e^{(A'\lambda)_j}-1\right]\,\geq\,0\]
(because all the $Q_j$ have nonnegative coefficients), which implies
$b'\lambda\geq0$. Hence, we retrieve that $b'\lambda\ge0$ 
whenever $A'\lambda\geq0$, which is to be expected since of course,
the existence of nonnegative 
integral solutions to (\ref{intro1}) implies the existence of
nonnegative real solutions.
\vspace{0.2cm}

%\noindent
{\bf Methodology}. We use counting techniques based on generating functions
as described by Barvinok and Pommersheim in \cite{barvinok} and by
Brion and Vergne in \cite{brion1,brion2}, to easily obtain 
a simple explicit expression of
the generating function (or, $\Z$-transform) $F:\C^m\to\C$ of the function 
$f:\Z^m\to\N$, $b\mapsto\,f(b)$, that counts the lattice points $x\in\N^n$ of the
convex polytope $\Omega:=\{x\in\R^n_+\,|\,Ax=b\}$.
Then $f$ is the inverse $\Z$-transform of $F$ and can be calculated
by a complex integral. Existence of a solution $x\in\N^n$ to (\ref{intro1})
is equivalent to showing that $f(b)\geq 1$, and by
a detailed analysis of this complex integral,
we prove that (\ref{intro2}) is a 
{\it necessary and sufficient} condition on $b$ for $f(b)\geq 1$.

\section{Notation and preliminary results}

For a vector $b\in\R^m$ and a matrix $A\in \R^{m\times n}$, denote by
$b'$ and $A'\in\R^{n\times m}$ their respective transpose. 
Denote by $e_m\in\R^m$ the vector with all entries equal to $1$.
Let $\R[x_1,\ldots,x_n]$ be the ring of real-valued polynomials in 
the variables $x_1,\ldots,x_n$. A polynomial 
$f\in\R[x_1,\ldots,x_n]$ is written
\[x\,\mapsto\,f(x)\,=\,\sum_{\alpha\in\N^n}f_\alpha x^\alpha\,=\,
\sum_{\alpha\in\N^n}f_\alpha x_1^{\alpha_1}\cdots x_n^{\alpha_n},\]
for {\it finitely many} real coefficients $\{f_\alpha\}$.

Given a matrix $A\in\Z^{m\times n}$, let $A_j\in\Z^m$ denote its $j$-th column
(equivalently, the $j$-th row of $A'$); then for every $z\in\C^m$, $z^{A_j}$ stands for 
\[z^{A_j}\,:=\,z_1^{A_{1j}}\cdots z_m^{A_{mj}}\,=\,
\e^{\langle A_j,\ln{z}\rangle}\,=\,\e ^{(A'\ln{z})_j}.\]
If $A_j\in\N^m$ then $z^{A_j}$ is a monomial of $\R[z_1,\ldots,z_m]$.

\subsection{Preliminary result}

Let $A\in\Z^{m\times n},\,b\in\Z^m$ and consider the system of linear equations
\begin{equation}
\label{2-1}
Ax\,=\,b;\quad\quad x\in\N^n,
\end{equation}
and its associated convex polyhedron 
\begin{equation}
\label{omega}
\Omega\,:=\,\{x\in\R^n\,\vert\,Ax=b;\:x\geq0\}.
\end{equation}
It is assumed that the {\it recession cone} $\{x\in\R^n\,\vert\,Ax=0;x\geq0\}$ 
of $\Omega$, reduces to the singleton $\{0\}$, so that $\Omega$ is compact
(equivalently, $\Omega$ is a convex polytope).

By a specialized version of a Farkas Lemma due to Carver,
(see e.g. Schrijver in \cite[(33), p. 95]{schrijver}),
this in turn implies that 
\begin{equation}
\label{u21}
\{\lambda\in\R^m\,\vert\,A'\lambda\,>\,0\}\,\neq\,\emptyset.
\end{equation}
Denote by $b\mapsto f(b)$ the function
$f:\Z^m\to\N$ that counts the nonnegative integral solutions $x\in\N^m$ of the 
system of linear equations (\ref{2-1}),
that is, the lattice points $x\in\N^n$ of $\Omega$. In view of (\ref{u21}), $f(b)$ is finite for all $b\in\Z^m$ because
$\Omega$ is compact.
Let $F:\C^m\,\to\,\C$ be the two-sided
$\Z$-transform of $f$, that is,
\begin{equation}
\label{2-2}
z\,\mapsto\,F(z)\,:=\,\sum_{u\in\Z^m}f(u)\,z^{-u}\,=\,
\sum_{u\in\Z^m}f(u)\,z_1^{-u_1}\cdots z_m^{-u_m}
\end{equation}
when the above series converges on some domain $D\subset\C^m$.
It turns out that $F(z)$ is well-defined on
\begin{equation}
\label{dom}
D\,:=\,\{z\in\C^m\,\vert\quad\vert z_1^{A_{1j}}\cdots z_m^{A_{mj}}\vert\,>\,1\quad 
j=1,\ldots,n\}.
\end{equation}

\begin{proposition}
\label{prop1}
Let $A\in\Z^{m\times n},b\in\Z^n$ and assume that (\ref{u21}) holds. Then :
\begin{equation}
\label{prop1-1}
F(z)\,=\,\frac{1}{\prod_{j=1}^n(1-z^{-A_j})}\,=\,
\frac{1}{\prod_{j=1}^n(1-z_1^{-A_{1j}}\cdots z_m^{-A_{mj}})}
\end{equation}
for all $z\in\Z^m$ that satisfy
\begin{equation}
\label{prop1-2}
\vert z^{A_j}\vert\,=\,\vert z_1^{A_{1j}}\cdots z_m^{A_{mj}}\vert\,>\,1\quad\quad j\,=\,1,\ldots,n.
\end{equation}
Moreover,
\begin{eqnarray}
\label{prop1-3}
f(b)&=&\frac{1}{(2\pi i)^m}\int_{\vert z_1\vert=\gamma_1}\cdots
\int_{\vert z_m\vert =\gamma_m}\frac{z^{b-e_m}}{\prod_{j=1}^n(1-z_1^{-A_{1j}}
\cdots z_m^{-A_{mj}})}\,dz\\
\nonumber
&=&\frac{1}{(2\pi i)^m}\int_{z\in\Gamma}\:\frac{z^{b-e_m}}{\prod_{j=1}^n(1-z^{-A_{j}})}
\,dz,
\end{eqnarray}
with $\Gamma\,:=\,\{z\in\C^m\,|\,\vert z_j\vert=\gamma_j\}$, and where
$\gamma\in\R^m_+$ is fixed and satisfies
\begin{equation}
\label{domain}
\gamma^{A_j}\,=\,\gamma_1^{A_{1j}}\cdots \gamma_m^{A_{mj}}\,>\,1\quad\quad j=1,\ldots,n.
\end{equation}
\end{proposition} 
\begin{proof}
The proof is a verbatim copy of that of Lasserre and Zeron in \cite{laszeron1}
where the linear system $Ax\leq b$ (instead of $Ax=b$) was considered, but for
the sake of completeness we reproduce it here.
Apply the definition (\ref{2-2}) of $F$ to obtain :
\begin{eqnarray*}
F(z)&=&\sum_{u\in\Z^m}z^{-u}
\left[\sum_{x\in\N^n,\;Ax=u}1\right]\,=\,
\sum_{x\in\N^n}\left[\sum_{u= Ax}z_1^{-u_1}
z_2^{-u_2}\cdots z_m^{-u_m}\right]\\
&=&\sum_{x\in\N^n}z_1^{-(Ax)_1}z_2^{-(Ax)_2}\cdots z_m^{-(Ax)_m},
\end{eqnarray*}
Now observe that 
$$z_1^{-(Ax)_1}z_2^{-(Ax)_2}\cdots z_m^{-(Ax)_m}\,=\,\prod_{k=1}^m
\left(z_1^{-A_{1k}}z_2^{-A_{2k}}\cdots z_m^{-A_{mk}}\right)^{x_k}\,=\,
\prod_{k=1}^m\left(z^{-A_{k}}\right)^{x_k}.$$
Hence, when (\ref{prop1-2}) holds we obtain
\begin{eqnarray*}
F(z)&=&\prod_{k=1}^n\,\sum_{x_k=0}^\infty\left( 
z^{-A_{k}}\right)^{x_k}\,=\,
\prod_{k=1}^n\left[1-z^{-A_{k}}\right]^{-1},
\end{eqnarray*}
which is (\ref{prop1-1}), and (\ref{prop1-3}) is obtained by a direct application of the inverse 
$\Z$-transform (see e.g. Conway in \cite{conway}). 
It remains to show that, indeed, the domain  defined in (\ref{prop1-2}) is not empty.
But this follows from (\ref{u21}). Indeed take $z_k:=\e ^{\lambda_k}$ for all
$k=1,\ldots ,m$, for any $\lambda$ that satisfies (\ref{u21}).
\end{proof}

\section{Main result}

Before proceeding to the general case $A\in\Z^{m\times n}$, we first consider
the case $A\in\N^{m\times n}$ where $A$ (and thus $b$) has only nonnegative entries.

\subsection{The case $A\in\N^{m\times n}$}
\label{subcase}

In this section $A\in\N^{m\times n}$ and  thus, necessarily $b\in\N^m$
(otherwise $\Omega=\emptyset$).
\begin{theorem}
\label{th2}
Let $A\in\N^{m\times n},b\in\N^m$. Then the following two statements (i) and
(ii) are equivalent :

{\rm (i)} The linear system $Ax=b$ has a solution $x\in\N^n$.

{\rm (ii)} The real-valued polynomial $z\,\mapsto\,z^b-1:=z_1^{b_1}\cdots z_m^{b_m}-1$ can be written
\begin{equation}
\label{th2-1}
z^b-1\,=\,\sum_{j=1}^nQ_j(z)(z^{A_j}-1)
%\,=\,\sum_{j=1}^nQ_j(z)(z_1^{A_{1j}}z_2^{A_{2j}}\cdots
%z_m^{A_{mj}}-1)
\end{equation}
for some real-valued polynomials $Q_j\in\R[z_1,\ldots,z_m]$, $j=1,\ldots,n$, all
of them with nonnegative coefficients.

In addition, the degree of the $Q_j$'s in (\ref{th2-1}) is bounded by
\begin{equation}
\label{b*}
b^*\,:=\,\displaystyle{\sum_{j=1}^mb_j-\min_k\sum_{j=1}^mA_{jk}}.
\end{equation}
\end{theorem}
For a proof see \S\ref{proofs}.

\subsection{Discussion}
(a) 
Let  $s(u):={m+u\choose u}$
the dimension of the vector space of
polynomials of degree $u$ in $m$ variables.
In view of Theorem \ref{th2}, and with $b^*$ as in (\ref{b*}), checking the existence of a solution 
$x\in\N^n$ to $Ax=b$ reduces to checking whether or not there exists a nonnegative solution $y$
to a system of linear equations with :
\begin{itemize}
\item $n\times s(b^*)$ variables, the nonnegative coefficients of the $Q_j$'s.

\item $s(b^*+\displaystyle{\max_k\sum_{j=1}^nA_{jk}})$ equations to identify the terms of same power
in both sides of (\ref{th2-1}).
\end{itemize}
This in turn reduces to solving a LP problem with $ns(b^*)$ variables and 
$s(b^*+\max_k\sum_jA_{jk})$ equality constraints. Observe that in view
of (\ref{th2-1}), this LP
has a matrix of constraints with only $0$ and $\pm1$ coefficients.

(b) In fact, from the proof of Theorem \ref{th2}, it
follows that one may even 
enforce the weights $Q_j$ in (\ref{th2-1}) to be polynomials in
$\Z[z_1,\ldots,z_m]$ (instead of $\R[z_1,\ldots,z_m]$) with
nonnegative coefficients (and even with coefficients in $\{0,1\}$).
However, (a)  above shows that the strength of Theorem \ref{th2} is
precisely to allow $Q_j\in\R[z_1,\ldots,z_m]$ as it permits to 
check feasibility by solving a (continuous) linear program.
Enforcing $Q_j\in\Z[z_1,\ldots,z_m]$ would result in an {\it integer}
program of size larger than that of the original problem.
%\end{remark}

(c) Theorem \ref{th2} reduces the issue of existence of a solution $x\in\N^n$ to
a particular {\it ideal membership problem}, that is, $Ax=b$ has a
solution $x\in\N^n$ if and only if the polynomial
$z^b-1$ belongs to the {\it binomial ideal}
$I=\langle z^{A_j}-1\rangle_{j=1,\ldots ,n}\subset\R[z_1,\ldots,z_m]$ 
{\it and} for some weights $Q_j$ all with {\it nonnegative coefficients}.

Interestingly, consider the ideal $J\subset\R[z_1,\ldots,z_m,y_1,\ldots,y_n]$ generated by
the binomials $z^{A_j}-y_j$, $j=1,\ldots,n$, and let $G$ be a 
Gr\"obner basis of $J$. Using the algebraic approach described by
Adams and Loustaunau in \cite[\S 2.8]{adams},
it is known that $Ax=b$ has a solution $x\in\N^n$ if and only if the monomial
$z^b$ is reduced (with respect to $G$) to some monomial $y^\alpha$, in which case
$\alpha\in\N^n$ is a feasible solution. Observe that this is not a
Farkas lemma as we do not know in advance $\alpha\in\N^n$ 
(we look for it!) to test whether
$z^b-y^\alpha\in J$. One has to apply Buchberger's algorithm to (i)
find a reduced Gr\"obner basis $G$ of $J$, and (ii) reduce $z^b$ with
respect to $G$ and check whether the final result is a monomial
$y^\alpha$. Moreover, note that the latter approach uses
polynomials in $n+m$ (primal) variables $y$ and (dual) variables $z$,
in contrast with the (only) $m$ dual variables $z$ in Theorem \ref{th2}.

\subsection{The general case}

In this section we consider the general case $A\in\Z^{m\times n}$ 
so that $A$ may have negative entries.
The above arguments cannot be repeated because of the occurence of negative powers.
However, let $\alpha\in\N^n,\beta\in\N$ be such that
\begin{equation}
\label{3-1}
\widehat{A}_{jk}\,:=\,
A_{jk}+\alpha_k\,\geq \,0;\quad \widehat{b}_j\,:=\,b_j+\beta\,\geq\,0;\quad k=1,\ldots,n;\,j=1,\ldots,m.
\end{equation}
Note that once $\alpha\in\N^n$ is fixed as in (\ref{3-1}), 
we can choose $\beta\in\N$ as large as desired.
Moreover, as $\Omega$ defined in (\ref{omega}) is compact, we have
\begin{equation}
\label{om1}
\max_{x\in\N^n}\,\{\sum_{j=1}^n\alpha_jx_j\,\vert\,Ax=b\}\,\leq\,
\max\,\{\sum_{j=1}^n\alpha_jx_j\,\vert\,x\in\Omega\}\,=:\,\rho ^*(\alpha)\,<\,\infty.
\end{equation}
Given $\alpha\in\N^n$, the scalar $\rho^*(\alpha)$ is easily calculated
by solving a LP problem. Next, choose $\rho^*(\alpha)\leq\beta\in\N$, and
let $\widehat{A}\in\N^{m\times n},\widehat{b}\in\N^{m}$ be as in (\ref{3-1}).
Then the existence of solutions $x\in\N^n$ to $Ax=b$ is equivalent to the existence
of solutions $(x,u)\in\N^n\times \N$ to the system of linear equations
\begin{equation}
\label{3-2}
\Q\:\left\{\begin{array}{rclcl}
\widehat{A}x&+&u\,e_m&=&\widehat{b}\\
\displaystyle{\sum_{j=1}^n\alpha_jx_j}&+&u&=&\beta.\end{array}\right.
\end{equation}
Indeed, if $Ax=b$ with $x\in\N^n$, then
\[Ax+e_m\sum_{j=1}^n\alpha_jx_j-e_m\sum_{j=1}^n\alpha_jx_j\,=\,b+(\beta-\beta)\,e_m,\]
or, equivalently,
\[\widehat{A}x +\left(\beta-\sum_{j=1}^n\alpha_jx_j\right)e_m\,=\,\widehat{b},\]
and thus, as $\beta\geq\rho^*(\alpha)\geq\sum_{j=1}^n\alpha_jx_j$ (cf. (\ref{om1})), 
we see that $(x,u)$ with $\beta-\sum_{j=1}^n\alpha_jx_j=:u\in\N$, is a solution of (\ref{3-2}).
Conversely, let $(x,u)\in\N^n\times\N$ be a solution of (\ref{3-2}).
Then, using the definitions of $\widehat{A}$ and $\widehat{b}$,
\[Ax+e_m\sum_{j=1}^n\alpha_jx_j +u\,e_m=b+\beta \,e_m;\quad\sum_{j=1}^n\alpha_jx_j+u\,=\,\beta,\]
so that $Ax=b$.
The system of linear equations (\ref{3-2}) can be put in the form
\begin{equation}
\label{equi}
B\left[\begin{array}{c}x\\u\end{array}\right]\,=\,
\left[\begin{array}{c}\widehat{b}\\ \beta\end{array}\right]
\mbox{ with}\quad B\,:=\,
\left[\begin{array}{rcl}\widehat{A}&|&e_m\\ 
-&&-\\
\alpha '&|&1\end{array}\right],
\end{equation}
and as $B\in\N^{(m+1)\times (n+1)}$, we are back to the case analyzed
in \S\ref{subcase}.

\begin{theorem}
\label{th3}
Let $A\in\Z^{m\times n},b\in\Z^m$ and assume that 
$\Omega$ defined in (\ref{omega}) is compact.
Let $\widehat{A}\in\N^{m\times n},\widehat{b}\in\N^m$, $\alpha\in\N^n$ and
$\beta\in\N$ be as in (\ref{3-1})  with $\beta\geq\rho^*(\alpha)$
(cf. (\ref{om1})). Then the following two statements (i) and (ii) are
equivalent :

{\rm (i)} The system of linear equations $Ax=b$ has a solution
$x\in\N^n$.

{\rm (ii)} The real-valued polynomial $z\mapsto 
z^b(zy)^\beta-1\in\R[z_1,\ldots,z_m,y]$ can be written
\begin{equation}
\label{th3-1}
z^b(zy)^\beta-1\,=\,Q_0(z,y)(zy-1)+\sum_{j=1}^nQ_j(z,y)(z^{A_j}(zy)^{\alpha_j}-1)
\end{equation}
for some real-valued polynomials $\{Q_j\}_{j=0}^n$ in
$\R[z_1,\ldots,z_m,y]$, all with nonnegative coefficients.

In addition, the degree of the $Q_j$'s in (\ref{th3-1}) is bounded by
\[(m+1)\beta +\sum_{j=1}^mb_j-\min\left[m+1,\,\min_{k=1,\ldots,n}\left[(m+1)\alpha_k+\sum_{j=1}^mA_{jk}\right]\right].\]

\end{theorem}
\begin{proof}
Apply Theorem \ref{th2} to the equivalent form (\ref{equi}) 
of the system $\Q$ in (\ref{3-2}) (since
$B\in\N^{(m+1)\times (n+1)}$ and $(\widehat{b},\beta)\in\N^{m+1}$),
and use the definition (\ref{3-1}) of $(\widehat{b},\beta)$ and
$\widehat{A}$.
\end{proof}

\section{Proof of Theorem \ref{th2}}
\label{proofs}

\begin{proof}
(ii) $\Rightarrow$ (i). Assume that $z^b-1$ can be written as in
(\ref{th2-1}) for some polynomials $\{Q_j\}$ with nonnegative
coefficients $\{Q_{j\alpha}\}$, that is,
\[Q_j(z)\,=\,\sum_{\alpha\in\N^m}Q_{j\alpha}z^\alpha\,=\,
\sum_{\alpha\in\N^m}Q_{j\alpha}z_1^{\alpha_1}\cdots z_m^{\alpha_m},\]
for finitely many nonzero (and nonnegative) coefficients $\{Q_{j\alpha}\}$.
By Proposition \ref{prop1},
the number $f(b)$ of nonnegative integral solutions $x\in\N^n$ to the equation $Ax=b$, is
given by
\[f(b)\,=\,\frac{1}{(2\pi i)^m}\int_{\vert z_1\vert=\gamma_1}\cdots
\int_{\vert z_m\vert=\gamma_m}\frac{z^{b-e_m}}
{\prod_{j=1}^n (1-z^{-A_{k}})}\,dz.\]
Writing $z^{b-e_m}$ as $z^{-e_m}(z^b-1+1)$ we obtain
\[f(b)\,=\,B_1+B_2,\]
with
\[B_1\,=\,\frac{1}{(2\pi i)^m}\int_{\vert z_1\vert=\gamma_1}\cdots
\int_{\vert z_m\vert=\gamma_m}\frac{z^{-e_m}}
{\prod_{k=1}^n (1-z^{-A_{k}})}\,dz,\]
and 
\begin{eqnarray*}
B_2&:=&\frac{1}{(2\pi i)^m}\int_{\vert z_1\vert=\gamma_1}\cdots
\int_{\vert z_m\vert=\gamma_m}\frac{z^{-e_m}(z^{b}-1)}
{\prod_{k=1}^n (1-z^{-A_{k}})}\,dz\\
&=&\sum_{j=1}^n\frac{1}{(2\pi i)^m}\int_{\vert z_1\vert=\gamma_1}\cdots
\int_{\vert z_m\vert=\gamma_m}\frac{z^{A_j-e_m}Q_j(z)}
{\prod_{k\neq j} (1-z^{-A_{k}})}\,dz\\
&=&\sum_{j=1}^n\sum_{\alpha\in\N^m}\frac{Q_{j\alpha}}{(2\pi  i)^m}
\int_{\vert z_1\vert=\gamma_1}\cdots
\int_{\vert z_m\vert=\gamma_m}\frac{z^{A_j+\alpha-e_m}}
{\prod_{k\neq j} (1-z^{-A_{k}})}\,dz.
\end{eqnarray*}
From (\ref{prop1-3}) in Proposition \ref{prop1} (with $b:=0$)
we recognize in $B_1$ the number of solutions $x\in\N^n$ to the 
linear system $Ax=0$, so that $B_1=1$.
Next, again from (\ref{prop1-3}) in Proposition \ref{prop1} 
(now with $b:=A_j+\alpha$), each term 
\[C_{j\alpha}\,:=\,\frac{Q_{j\alpha}}{(2\pi  i)^m}
\int_{\vert z_1\vert=\gamma_1}\cdots
\int_{\vert z_m\vert=\gamma_m}\frac{z^{A_j+\alpha-e_m}}
{\prod_{k\neq j} (1-z^{-A_{k}})}\,dz,\]
is equal to
\[Q_{j\alpha}\times \mbox{ the number of integral solutions }x\in\N^{n-1}\]
of the linear system
$\widehat{A}^{(j)}x=A_j+\alpha$, where
$\widehat{A}^{(j)}$ is the matrix in $\N^{m\times (n-1)}$ obtained
from $A$ by deleting its $j$-th column $A_j$. As by hypothesis, each
$Q_{j\alpha}$ is nonnegative, it follows that 
\[B_2\,=\,\sum_{j=1}^n\sum_{\alpha\in\N^m} C_{j\alpha}\,\geq\,0,\]
and so $f(b)=B_1+B_2\geq 1$. In other words,  the sytem $Ax=b$ has at least one solution $x\in\N^n$.

(i) $\Rightarrow$ (ii). Let $x\in\N^n$ be a solution of $Ax=b$, and write
\[z^{b}-1\,=\,z^{A_1x_1}-1+z^{A_1x_1}(z^{A_2x_2}-1)+\cdots+z^{\sum_{j=1}^{n-1}A_jx_j}(z^{A_nx_n}-1),\]
and
\[z^{A_jx_j}-1\,=\,(z^{A_j}-1)\left[1+z^{A_j}+\cdots+z^{A_j(x_j-1)}\right]\quad
j=1,\ldots,n,\]
to obtain (\ref{th2-1}) with
\[z\,\mapsto\,Q_j(z)\,:=\,z^{\sum_{k=1}^{j-1}A_kx_k}\left[1+z^{A_j}+\cdots+z^{A_j(x_j-1)}\right],\quad
j=2,\ldots,n,\]
and
\[z\,\mapsto\,Q_1(z)\,:=\,1+z^{A_1}+\cdots+z^{A_1(x_1-1)}.\]
We immediately see that each $Q_j$ has all its coefficients 
nonnegative (and even in $\{0,1\}$). 

Finally, the bound on the degree follows immediately from the
expression of the $Q_j$'s in the proof of (i) $\Rightarrow$ (ii).
\end{proof}

\end{document}